\documentclass[12pt]{article}
\usepackage{amssymb,amsmath,amscd}
\title{On Quasihomogeneous Curves}
\date{\mbox{}}
\begin{document}
\author{ By Abdallah Assi and Avinash Sathaye}

\thanks{
Universit\'e d'Angers, Math\'ematiques, 49045 Angers cedex 01, France;\\
\hspace*{4mm}
e-mail: assi@univ-angers.fr\\
\hspace*{4mm}
University of Kentucky, Lexington KY 40506, U.S.A.;\\
\hspace*{4mm}
email:sathaye@uky.edu
During the development of this work, Sathaye visited Universit\'e d'Angers; 
he would like to thank that institution for hospitality and support.}

\newcommand{\numero}{\refstepcounter{teorema}
\paragraph{\bf\theteorema}}

\newtheorem{teorema}{Theorem}[section]

\renewcommand{\theequation}{\theteorema}
\newenvironment{formule}{\refstepcounter{teorema}
\begin{eqnarray}}{\end{eqnarray}}
\def\nz{\hbox{\text{0}{ \raise 3.5pt \hbox{\kern -1pt
\vrule width 8pt height 0.5 pt }\kern 3pt}}}


\newtheorem{p}[teorema]{Theorem}
\newtheorem{proposicion}[teorema]{Proposition}
\newtheorem{lema}[teorema]{Lemma}
\newtheorem{definicion}[teorema]{Definition}
\newtheorem{consecuencia}[teorema]{Consecuencia}
\newtheorem{corolario}[teorema]{Corollary}
\newtheorem{ejercicio}[teorema]{Ejercicio}
\newtheorem{ejemplo}[teorema]{Ejemplo}
\newtheorem{cuestion}[teorema]{Cuesti\'on}
\newtheorem{nota}[teorema]{Remark}
\newtheorem{exemple}[teorema]{Example}
\newtheorem{theorema}[teorema]{Theorem}
\newenvironment{demostracion}[1]{\paragraph{\sl Proof#1}}{}
\newenvironment{notacion}[1]{\paragraph{\sl Notaci\'on#1}}{}
\newenvironment{remarque}[1]{\paragraph{\sl Remark#1}}{}
\newenvironment{descit}[1]{\begin{quote}{\em #1\/}:}{\end{quote}}
\newenvironment{proof}[1]{\paragraph{\sl Proof#1}}{\qed}
\newcommand{\jitem}{\stepcounter{enumi}\item[{\rm
\theenumi}.]}
\def\un{\underline}
\def\ov{\overline}
\def\Deltabar{\ov{\Delta}}
\def\KKn{{\bf K}^n}
\def\KKFX{{\KK[[\un{x}]]}}
\def\KKCX{{\KK\{\un{x}\}}}
\def\KKX{{\KK[\un{x}]}}
\def\KKY{{\KK[\un{y}]}}
\def\DDFK{{\widehat{\cD}_{n}(\KK)}}
\def\DDF{{\widehat{\cD}}}
\def\AnK{{A_n(K)}}
\def\PNK{{\PP_n(\KK)}}

\def\DK{\cD_{n}(\KK)}
\newcommand{\DD}{{\bf D}}
\newcommand{\bfm}{{\bf m}}
\newcommand{\bfp}{{\bf p}}
\newcommand{\bfq}{{\bf q}}
\newcommand{\CC}{{\bf C}}
\newcommand{\RR}{{\bf R}}
\newcommand{\KK}{{\bf K}}
\newcommand{\QQ}{{\bf Q}}
\newcommand{\ZZ}{{\bf Z}}
\newcommand{\NN}{{\bf N}}
\newcommand{\PP}{{\bf P}}
\newcommand{\HH}{{\bf H}}
\newcommand{\FF}{{\bf F}}
\newcommand{\LL}{{\bf L}}
\renewcommand{\SS}{{\bf S}}
\newcommand{\Dp}{{\mathop{\bf D}}}
\newcommand{\Cp}{{\mathop{\bf C}}}
\newcommand{\Rp}{{\mathop{\bf R}}}
\newcommand{\Qp}{{\mathop{\bf Q}}}
\newcommand{\Zp}{{\mathop{\bf Z}}}
\newcommand{\Np}{{\mathop{\bf N}}}
\newcommand{\extL}{{\prec_L}}
\newcommand{\extegalL}{{\preceq_L}}

\newcommand{\qed}{\mbox{$\Box$}}

\newcommand{\expp}{{\rm exp}}
\newcommand{\Expp}{{\rm Exp}}
\def\New{\mbox{\it Newton}}
\newcommand{\resp}{{\em resp.}}
\newcommand{\ppcm}{{\em ppcm}}
\newcommand{\see}{{\em voir}}
\newcommand{\conf}{{\em cf.}}
\newcommand{\ie}{{\em i.e.}}
\newcommand{\loccit}{{\em loc. cit.}}

\newcommand{\val}{{\rm val}}
\newcommand{\ord}{{\rm ord}}
\newcommand{\gr}{{\rm gr}}



\begin{abstract} A hypersurface is said to be quasihomogeneous if in suitable coordinates with
assigned weights, its equation becomes weighted homogeneous in its
variables. For an irreducible quasihomogeneous plane curve, the equation 
necessarily becomes a two term equation of the form $aY^n+bX^m$ where $n,m$ are
necessarily coprime. Zariski, in a short paper, established a criterion 
for an algebroid
curve to be quasihomogeneous \cite{Z1} and a celebrated theorem of 
Lin and Zaidenberg gives a global criterion for quasihomogeneity \cite{LZ}.
The Lin-Zaidenberg theorem does not have a simple proof, despite having
three different proofs using function theory, topology and algebraic
surface theory respectively. We give here a global version of the
Zariski result. As a consequence we give a proof of a slightly weaker
version of the Lin-Zaidenberg Theorem, namely that a rational curve with 
one place at infinity is unibranch and locally quasihomogeneous if and only if it is 
globally quasihomogeneous, provided the ground field is algebraically
closed of characteristic zero.
Our method of proof leads to some interesting questions about the
change in the module of differentials when we go to the integral
closure.

\end{abstract}
\maketitle

\medskip

\section{Introduction and notation}
 
Let ${\bf k}$ be an algebraically closed field of characteristic 
$0$ and let $A=k^{[2]} = k[X,Y]$ be a polynomial ring in two variables over $k$. 
An irreducible polynomial $f=f(X,Y)\in A$ defines an irreducible plane curve. 
We shall denote its coordinate ring by $R=A/(f)$ and the canonical
homomorphism $A\longrightarrow R$ shall be denoted by $\phi_f$.
In general we write $\phi_f(X)=x,\phi_f(Y)=y$. 
For any $h=h(X,Y)\in A$ we shall denote by $J_f(h)=J_f(h(X,Y))$ its jacobian with
$f$, namely $h_Xf_Y-h_Yf_X$. 
The ideal in $A$ generated by all such jacobians
is denoted by $J_f=(f_X,f_Y)A$. We also need the extended ideal
$Jac_f = (f,J_f)A$. 
Note that for any $u=g/h\in qt(A)$ the quotient field of $A$, we can
easily extend the definition to calculate $J_f(u)$, by the usual rules
of derivatives.
We shall denote by $\overline{J}_f $ the ideal $\phi_f(J_f)\subset R$. 
\medskip

\noindent
{\bf We may often drop the reference to $f$ if it is clear during
 a discussion and
we may drop explicit reference to $\phi_f$ as follows.}
We may write $J_f(h(x,y))$ or simply $J(h(x,y))$ to denote the image
$\phi_f(h_X(X,Y)f_Y-h_Y(X,Y)f_X)$. 
If we consider $u(X,Y)=g(X,Y)/h(X,Y)$ where $\phi_f(h(X,Y)) = h(x,y)\neq
0$, we have a well defined element of $k(f)$: 
$\phi_f(J_f(u)) = \phi_f(J_f(g))/\phi_f(h) -
\phi_f(gJ_f(h))/\phi_f(h)^2$. We have thus defined $J_f(u(x,y)$ 
for any $u(x,y)\in k(f)$.

\noindent
We introduce the following additional notations for future use.
Define $\chi(f;Q)$ to be the number of branches of $f$ at $Q$,
and let $\overline\chi(f;Q)= \chi(f;Q)-1$. 
Let $\overline{R}$ denote the integral closure 
of $R$ in its quotient field,
i.e. the function field $k(f)$ of the curve $f$. We define the
conductor ideal 
$\mathfrak{C}(\overline{R},R)=\{u \in R ~|~ u\overline{R} \subset R\}$ and
may shorten the notation to $\mathfrak{C}_f$ or simply $\mathfrak{C}$ for
convenience.

Let also
$\mathfrak{C}(f;Q)$ be the conductor ideal of 
$R_{Q'}$ in its integral closure (denoted as $\overline{R}_{Q'}$),
where $Q'=\phi_f(Q)$. Let  $C(f;Q)$ be the length 
$l(\overline{R}_{Q'}/ \mathfrak{C}(f;Q))$ as $R$-modules.

We define $C(f)$ to be the sum $\sum_Q C(f;Q)$ over all points $Q$ such
that $f$ passes thru $Q$. 

As is well known, $\mathfrak{C}_f$ is  an ideal in $R$ as well as 
$\overline{R}$ and has the property 
$l(\overline{R}/\mathfrak{C}) = 2l(R/\mathfrak{C})$.
Further this is clearly equal to $C(f)$.

\medskip

\noindent
We have the usual module $\Omega(R,k)=\Omega(R)$ of 1-differentials
generated by $dx,dy$ and subject to the relation $f_xdx+f_ydy =0$. 
If  $h(x,y)\in k(f)$, such that $dh \neq 0$, then the differential
$dh(x,y)/J_f(h(x,y))$ is easily seen to be independent of the choice of $h$ when
viewed as a member of
$\Omega(k(f),k)$. We will denote it by $\omega$ and call it the
fundamental differential of $f$.
In our zero characteristic case, the condition $dh \neq 0$ can be
replaced by ``$h$ is a non constant''.

Given non constant $u\in \overline{R}$, we have $du \neq 0$.
For any valuation $w$ of the function field $k(f)$ over $k$, 
recall that the value $w(du)$ of the differential 
$du$ is defined as follows. Choose a parameter $t$ for the valuation $w$
which is a separating transcendence 
(separability being automatic in our characteristic
zero). We have $du/J_f(u) = dt/J_f(t)$ and we define
$w(du)=w(J_f(u)/J_f(t))$. It is then known that this is independent of
the choice of the parameter $t$. In short, this can be also described as 
$w(du)=w(h)$ where $h$ is chosen such that $du=hdt$ for some
(separating) uniformizing parameter $t$ for $w$.

\section{\bf Generalized Berger formula}

\begin{lema} 
Let $f\in A$ be a non constant irreducible polynomial and using the
above notation, let
$$
J^* =( \{ J_f(u) ~|~ u \in \overline{R} \}) \overline{R}
$$

Then we have:

$$J^*=\mathfrak{C}=\mathfrak{C}_f.$$
\end{lema}

\begin{demostracion}{.} 

\medskip

Let $w$ be any valuation at finite distance (i.e. the valuation ring of
$w$ contains $R$.
By the Dedekind formula for the conductor and differential \cite{AS}(28.15.1) 
we know that 
$w(J_f(u))=w(\mathfrak{C})+w(du)$. 
 
Since $u \in \overline{R}$ we get $w(du) \geq 0$ and 
hence $w(J_f(u))\geq w(\mathfrak{C})$.
Thus $w(J^*) \geq w(\mathfrak{C})$. Further, if we choose $u$ to be itself a
uniformizing parameter for $w$ in $\overline{R}$, then
$w(J_f(u))=w(\mathfrak{C})$ and thus the minimum of values of elements of
$J^*$ at $w$ is actually equal to $w(\mathfrak{C})$.
Thus $w(J^*) = w(\mathfrak{C})$ for all valuations dominating the normal domain 
$\overline{R}$ and  hence $J^*=\mathfrak{C}$.

\end{demostracion}

\begin{corolario}  With the same setup as in Lemma 2.1, we have:
${\mathfrak C}\omega = \Omega(\overline{R},k)$, where we 
recall that $\omega$ is the fundamental differential of $f$ as described
above.
\end{corolario}

\begin{demostracion}{.} 
Note that $\Omega(\overline{R},k)$ consists of sums of elements $adu$
with $a,u\in \overline{R}$. We have $adu = aJ_f(u)\omega \in J^*\omega$. Thus
$\Omega(\overline{R},k) \subset J^*\omega$.
Similarly $J^*\omega$ consists of sums of elements $aJ_f(u)\omega$
for $a,u \in \overline{R}$. But $aJ_f(u)\omega =
aJ_f(u)(du/J_f(u))=adu \in \Omega(\overline{R},k)$. 
Thus $J^*\omega \subset \Omega(\overline{R},k)$.

\end{demostracion}

\noindent By a similar argument we get the following:

\begin{lema} With the same setup as in Lemma 2.1, we have:
$ \overline{J}_f\omega = \Omega(R,k)$.
\end{lema}

\begin{demostracion}{.} 
Just imitate the argument of the above corollary.

\end{demostracion}

\noindent We can now state the following:

\begin{proposicion} We the notations as explained above, we have:
$$
l(A /{Jac_f})=l(R/ {\overline{J}_f})=
l(R/\mathfrak{C})+l(\mathfrak{C}/\overline{J}_f)=
l( \overline{R} / R)+l(\Omega(\overline{R},k)/ \Omega(R,k))
$$
\end{proposicion}

\begin{demostracion}{.} The first equality is obvious and the second is clear since
$\overline{J}_f$ is contained in $\mathfrak{C}$.

The multiplication by the fundamental differential $\omega$
clearly gives an injective homomorphism of $\mathfrak{C}$ into
$\Omega(\overline{R},k)$  and we get:
$$
l(\mathfrak{C}/\overline{J}_f)=l(\mathfrak{C}\omega / \overline{J}_f \omega)
$$

We use Corollary 2.2. and Lemma 2.3 to get our quantity equal to
$l(R/\mathfrak{C})+l(\Omega(\overline{R},k)/ \Omega(R,k))$. 
Finally, we use the well known fact that 
$ l(R/\mathfrak{C})= l(\overline{R} / R)$ to finish the proof.
\end{demostracion}

\begin{nota}{\rm We now connect the above calculations with 
Zariski's concept of the torsion module of differentials.
Let $T(R)$ be the torsion module of (the $R$-module) 
$\Omega(R,k)$, i.e., the submodule consisting of the 
elements of $\Omega(R,k)$ which have a non vanishing annihilator.
Explicitly, we see that: 
$$T(R) = \{ adx+bdy ~|~ a,b \in R, uadx+ubdy = 0 
\mbox { for some nonzero } u \in R\}.$$ 
\medskip

\noindent 
Zariski showed that in case $R$ is the local ring of an analytically
irreducible plane curve, the length of the torsion module satisfies:
$$l(T(R)) = l(R/\overline{J}_f).$$ 
where $f$ is the local equation.
Zariski's arguments remain valid for our global case as well and we get
that this result continues to hold.
\medskip

\noindent 
Zariski further used the Berger formula giving a local result similar to our
global statement of Proposition 2.4.
\medskip

\noindent 
Finally, Zariski noticed that the term $l(\Omega(\overline{R},k)/
\Omega(R,k))$ depends on the structure of the equation and is at most
equal to the adjacent term ( the length of the integral closure over
the coordinate ring or equivalently the length of the conductor ideal in the
coordinate ring). He showed that it reaches this maximum value 
if and only if under a suitable local analytic 
change of coordinates the curve becomes $y^m-x^n$ for coprime natural 
numbers $m,n$. 

The result depends on the analytic irreducibility, i.e. the property of
having only one place at the point. The main theme of our paper is to
globalize this result by showing that it continues to hold under the
assumption that the global curve has only one place at infinity. 

Saito has generalized the Zariski result in a different direction. 
His result states that a local hypersurface with an isolated singularity 
is (locally) quasihomogeneous if and only if its equation belongs to the
jacobian ideal. \cite{SA}} Unfortunately, his function theory proof has
not yet been transformed to an algebraic form.
\end{nota}

\section{\bf Algebraic invariants associated with the curve $f$}

\medskip

\noindent Let the notations be as in Section 2. In the 
following, we shall introduce invariants associated with 
the curve $f$. Let $Q \in j\mbox{-}Spec(A)$ (i.e. $Q$ is a maximal ideal 
in $A$. We also
describe it as a point in the plane).  
We recall the {\bf intersection multiplicity } of $f$ and $g$ at $Q$,
denoted by int$(f,g;Q)$, to be the length
of the $A_Q$-module $A_Q/(f,g)A_Q$. Note that this length is finite
exactly when $f,g$ don't have a common factor passing thru $Q$, or
equivalently, $(f,g)A_Q$ is primary for $QA_Q$. Also, it is positive,
if and only if both the curves pass thru the point.
We recall the total intersection multiplicity of $f$ with $g$ in $A$, denoted
by int$(f,g;A)= \sum \mbox{ int}(f,g;Q)$ taken over
all the points in the plane.

\noindent
{\bf Milnor numbers}: Let $f \in A$ be non constant.
We define
$$\mu(f;Q) = l(A_Q/J_fA_Q)=\mbox{ int}(f_X,f_Y;Q). $$
and call it the (local) Milnor number of $f$ at $Q$.
Note that this is finite if and only if $f$ is reduced, i.e.
does not have a multiple factor thru $Q$.
We also define $\mu(f)$ the affine Milnor number of $f$
as the sum of the local Milnor numbers over all points $Q$ such that
$f$ passes thru $Q$, i.e. $f\in
Q$. Note that this is finite if and only if $f$ has no multiple factors.
We define the Milnor number of the pencil $(f-a)_{a\in k}$
$\mu(f;A)$ by $\mu(f;A)=\sum_{a\in k}\mu(f-a)$.
Clearly, $\mu(f;A)= l(A/J_f)$ and $\mu(f;A) = \mu(f-a;A)$ for all $a\in
k$.
Note that this is finite if and only if all members of the pencil $\{f-a\}$
are reduced. 
\medskip

\noindent {\bf Tjurina numbers}: Let $f \in A$ be irreducible (non
constant).
We define
$$\nu(f;Q)=l(A_Q/Jac_f A_Q)\simeq 
l(R_{\phi_f(Q)}/(\overline{J}_f)R_{\phi_f(Q)})$$
and call it the local Tjurina number.
We also define $\nu(f)$ the affine Tjurina number 
of $f$ as the sum over all points $Q$ such that $f$ passes thru $Q$.
Clearly, $\nu(f) = l(A/Jac_f A) =  l(R/\overline{J}_fR)$.
The Tjurina numbers can be also defined for reduced $f$ and will still 
be finite.\footnote{We could also define a Tjurina number for a pencil, 
but so far no useful
consequences of such a definition are known, hence we refrain from
making such a definition. For the same reason, defect numbers below 
are also not defined for pencils.}

\noindent
In view of Zariski's work and our observations above, the Tjurina 
numbers could also be defined as lengths of the
corresponding torsion modules of the module of differentials (locally as
well as globally).

{\bf Defect numbers}: Let $f \in A$ be irreducible (non constant).
Since $J_f \subset Jac_f$, we define
$$\delta(f;Q) = \mu(f;Q)- \nu(f;Q) = l(J_fA_Q/Jac_fA_Q) $$ 
and call it the local defect number of $f$.  
We define $\delta(f)$ the affine defect of $f$ as the sum of local
defects $\delta(f;Q)$ taken over all the points $Q$ such that $f$ passes 
thru $Q$.

Thus: $\mu(f)=\delta(f)+\nu(f)$, i.e. Milnor=Tjurina+defect.

{\bf Zariski numbers}: Let $f \in A$ be non constant.
We define 
$Z(f;Q)= l(\mathfrak{C}(f;Q)/\overline{J}_fR_{Q'})$  where $Q'=\phi_f(Q)$.
We call it the local Zariski number of $f$ at $Q$.
We also define the affine Zariski number $Z(f)$ of $f$, 
as the sum of local Zariski numbers $Z(f;Q)$ taken over all the points
$Q$ such that $f$ passes thru $Q$. Note that then
$Z(f))=l(\mathfrak{C}_f/\overline{J}_f)$.

We now observe some useful relations between these numbers.
\begin{lema} 
Suppose $f\in A$ is irreducible (non constant) and passes thru a 
point $Q$. Then 
$$\mu(f;Q)+\overline{\chi}(f;Q)=C(f;Q).$$

\end{lema}

\begin{demostracion}{.}
The result is well known. We indicate an outline. 

First choose suitable local coordinates so that $Q$ is the maximal
ideal generated by $X,Y$ and that $X$ is not tangential to $f$ at $Q$, i.e.
at every valuation of $f$ centered at $Q$, the image $x = \phi_f(X)$ has 
minimal value in the ideal $(x,y)$. From the Dedekind formula we deduce that 

$$
C(f;Q) = int(f,f_Y;Q)-int(f,X;Q)+\chi(f;Q)\leqno{E1}.
$$
Next, by considering the coordinate ring of $f_Y$ and noticing that
modulo any irreducible factor of $f_Y$ we have $f_x=\frac{df}{dx}$, we
deduce that
$$ int(f,f_Y;Q) = int(f_X,f_Y;Q)+int(X,f_Y;Q)\leqno{E2}.$$
By nontangentiality of $X$, we further get that 
$$int(X,f_Y;Q) = int(X,f;Q)-1\leqno{E3}.$$

The result follows by combining these equations.
\end{demostracion}

\begin{lema} If $f\in A$ is irreducible (non constant) and passes thru 
a point $Q$ then we have
$$\mu(f;Q) = \nu(f;Q)+ \delta(f;Q) = C(f;Q)/2+Z(f;Q)+\delta(f;Q).$$

\end{lema}

\begin{demostracion}{.}
The first equality is by definition, while the second follows the
obvious fact that $$\nu(f;Q) =
l(R_{\phi_f(Q)}/\mathfrak{C}(f;Q))+Z(f;Q).$$ and the fact that the first
term is half the length $C(f;Q)$ of the conductor ideal in the integral
closure.
\end{demostracion}

\begin{lema} 
If $f\in A$ is unibranch and passes thru $Q$ then we have
$\mu(f;Q)=C(f;Q)$ and further
$$Z(f;Q)+\delta(f;Q)=C(f;Q)/2.$$
In particular, the maximum possible value for the local Zariski number
is $C(f;Q)/2$ and it is attained exactly when $\delta(f;Q)=0$. 

\end{lema}

\begin{demostracion}{.}
By Lemma 3.1, we get that $\mu(f;Q)=C(f;Q)$ and now Lemma 3.2 gives the
result.
\end{demostracion}

\begin{lema} 
If $f\in A$ is an irreducible (non constant) curve.
then we have the following:

$$C(f)=\mu(f)+\overline{\chi}(f) = C(f)/2+Z(f)+\delta(f)+ \overline{\chi}(f).$$
In particular,
$$C(f)/2-Z(f)=\delta(f)+ \overline{\chi}(f).$$

\end{lema}

\begin{demostracion}{.}
This is simply obtained by adding up the results from Lemmas 3.1 and
3.2 over all points of the curve $f$.
\end{demostracion}

\begin{nota}
{\rm 
Zariski's result can now be described as a necessary and sufficient
condition for maximality of $Z(f;Q)$ or vanishing of the defect.
To describe this, let us change coordinates so that $Q=(X,Y)$ and
consider $f\in A \subset k[[X,Y]]$. We say that $f$ is formally
quasihomogeneous at $Q$ if there is a change of variables
$k[[X,Y]]=k[[X',Y']]$ such that after giving some positive weights to
$X',Y'$ a local equation becomes a sum of equal weight monomials in $X',Y'$.
It is easy to see that by irreducibility of $f$ and a further change of
coordinates if necessary, we may assume that
$f=u((X')^m-(Y')^n)$ where $u$ is a unit in
the power series ring and $m,n$ are (necessarily comprime) positive integers.
We can now restate

\noindent {\bf Zariski Theorem} $\delta(f;Q)=\overline{\chi}(f;Q)=0$ if and only if $f$
locally irreducible and  formally quasihomogeneous at $Q.$



We shall be imitating Zariski's original proof in a global setting.
}
\end{nota}

\section{One place curves}

\medskip

\noindent Let the notations be as in Section 1. The curve $f$ is 
said to have one place at infinity if there is only one
valuation of its function field $k(f)$ over $k$ which does not contain
the coordinate ring $R$. Several important properties of such curves are
well known from the Abhyankar-Moh theory. We summarize them briefly.

\subsection{Setup}
Suppose that $f(X,Y)\not \in k[X]$ has one place at infinity (or equivalently
$f(X,Y) \in k((X))[Y]$ is irreducible).
Then it is essentially monic in all coordinate systems and can be written as

$$
f(X,Y)=\nz Y^n+a_1(X)Y^{n-1}+\cdots+a_n(X) \in k[X,Y]
$$

\noindent where $\nz$ is the well known "Abhyankar non zero", a symbol denoting
some non zero constant in $k$. Be aware that this symbol can denote different
quantities even in the same expression.
The highest $X$-degree term of $f$ can
then be described as $\nz X^m$.
We can and usually do assume that $f$ is actually monic in $Y$.

Let $v$ denote the unique valuation at infinity and denote by
$\Gamma(f)$ the semigroup consisting of all $v(h)$ as $h$ varies over non zero elements
of $R$. Note that $\Gamma(f)$ is a subset of $-N$ the set of negatives
of natural numbers. The structure and generation of this semigroup is
the center piece of the Abhyankar-Moh theory.

\subsection{Characteristic quantities}
There is a well defined sequence of numbers $r_i\in \Gamma(f)$ for $0
\leq i \leq h$ where
$$r_0=v(x)=-n, r_1=v(y)=-\deg_X(f(X,Y)).$$
Denoting the GCD$(r_0,\cdots,r_{i-1})$ by $d_i$ we get that $d_1=n \geq
d_2 > \cdots d_h >d_{h+1}=1$. Further these $d_i$ successively divide
the earlier ones and we define $n_i = d_i/d_{i+1}$ for $i=1,\cdots h$.
We additionally make {\bf  two special definitions}
$d_0=m$ and $n_0=d_0/d_2$.

There is a natural sequence of (essentially ) monic polynomials
$$g_0=X, g_1=Y, g_2(X,Y),\cdots,g_h(X,Y) $$
such that the $Y$-degrees of $g_i$ are $d_1/d_i=n/d_i$ for $i=1,\cdots h$.
Besides this degree condition, we also have that $v(g_i(x,y))=r_i$
for $i=0,\cdots h$.

Having made $f$ monic in $Y$, the $g_i$ may be chosen to be the
approximate $n/d_i$-th roots of $f$, which means, the unique polynomial of
$Y$-degree $n/d_i$ such that $\deg_Y(f-g_i^{d_i})<n-n/d_i$.  In general,
the chosen polynomials are called pseudo-approximate roots and are said
to form an associated $g$-sequence.

\subsection{Admissible expansions}
The polynomials give a standard basis for $R$ over $k$ in the
following sense.

A sequence $a=(a_0,\cdots,a_h)$ of integers is
said to be admissible if it satisfies the following two conditions
of admissibility:

{\bf Higher level condition:}\hspace*{0.3in} 
$n_i >a_i \geq 0$ for $i=1,\cdots h$ and

{\bf Level zero condition:} \hspace*{0.4in}$a_0 \geq 0$ 
Then by a standard $g$-monomial we mean an expression
$g^a= g_0^{a_0}\cdots g_h^{a_h}$, where $a=(a_0,\cdots,a_h)$ is
admissible.

The main theorem of the Abhyankar-Moh theory is that the set of standard
monomials form a $k$-basis for $R$ and moreover, the values
$v(g^a)=a_0r_0+\cdots +a_hr_h$ are distinct for distinct admissible $a$.

In particular
$$\Gamma(f)=\{ a_0r_0+\cdots +a_hr_h \mbox{ where }
(a_0,a_1,\cdots,a_h) \mbox{ is admissible } \}.$$

It can be shown by simple properties of integers that every
integer $w$ has a unique expression $w=\sum_0^h (b_ir_i)$, where the sequence
$b=(b_0,\cdots,b_h)$ satisfies higher level admissibility and $w \in
\Gamma(f)$ if and only if $b_0 \geq 0$ i.e. the zero level admissibility
holds as well.

\subsection{Properties of the Semigroup of values}
Define an integer:
$$\Theta=\Theta(f) = -r_0 + \sum_1^h (n_i-1)r_i.$$
As a consequence of the above description, easy calculations show that
given integers $p,q$ with $p+q=\Theta$ we have that exactly one of
$p,q \in \Gamma(f)$.

In particular, since $\Gamma(f)$ has no positive integers, all integers
less than $\Theta$ are in $\Gamma(f)$ while $\Theta\not \in \Gamma(f)$.
The integer $-1+\Theta$ is called the conductor of the semigroup
$\Gamma(f)$ and is necessarily an even number, since exactly half the
integers between $0$ and it are not in the semigroup. We shall denote it
by $C(\Gamma(f))$ and we have the formula:
$$ C(\Gamma(f)) = -1-r_0 + \sum_1^h (n_i-1)r_i .$$

\subsection{A conductor formula} 
Let us observe that this conductor of the semigroup is connected
with the length of the conductor ideal of the coordinate ring of the
curve by the formula:
$$C(\Gamma(f)) = -C(f)- 2P_g(f)$$
where $P_g(f)$ denotes the geometric genus of the function field $k(f)$.
Among its various equivalent definitions, let us choose the one which
says that it is the integer which equals $1$ plus half the degree of
the divisor of any nonzero differential in the function field. 
Explicitly, given any nonzero differential $adb$ we know that 
$2P_g(f)-2= \sum_w w(adb)$  where the sum is taken over all
valuations of the function field $k(f)$ over $k$ and the values $w(adb)$
are as described in section 1.

For our affine curve, we
can conveniently choose the fundamental differential $\omega = dx/f_y$ as 
described in section 2. 
The sum of values for valuations at finite distance is
easily seen to be $-C(f)$ by the Dedekind formula and hence the order
of  $\omega$ at the valuation $v$ at infinity is seen to be
$C(f)+2P_g(f)-2$.

On the other hand, we can compute it explicitly thus:
Since $v(x)=r_0$ we have $v(dx) = r_0-1$ and $v(f_y)$ is well known 
by the Abhyankar-Moh theory to be equal to $\sum_{1}^h (n_i-1) r_i$. 
\footnote{ This is worked out by starting with an NP expansion 
as described later in this section and explicitly determining 
the $t$-order of
the resultant $\mbox{Resultant}(f,f_Y,Y)$ interpreted as a product
$\prod (\eta(t)-\eta(\zeta t))$, where $\zeta$ varies over non identity
$n$-th roots of $1$.
}
Thus we get 
$$v(<\omega)>) = -1+r_0-\sum_{1}^h (n_i-1) r_i = C(f)+2P_g(f)-2.$$
Thus we have as claimed:
$$-C(f)-2P_g(f) = -1 - r_0 + \sum_{1}^h (n_i-1)r_i =
C(\Gamma(f)).$$
 
\subsection{\bf A formula for the Milnor number of the pencil.}
We recall that the Milnor number for the pencil $(f-a)_{a\in k}$ is
defined as $\mu(f;A) = l(A/J_f)$. By an argument similar to our proof in 
Lemma 3.1 we can argue that 
$$int(f_X,f_Y;A) = int(f,f_Y;A) - int(X,f_Y;A) =-\sum_{1}^h (n_i-1)r_i+ 1+r_0.$$
For proof, we note that the left hand side can be interpreted as the
negative of the sum of values at valuations at infinity for $f_Y$, with
suitable adjustment made for reducible $f_Y$. The fact that
$int(f-a,f_Y;A)$ is independent of $a$, being purely in terms of the
characteristic terms common to all translates $f-a$ and the fact that
$f$ is necessarily essentially monic in $Y$ guarantee that $f,x$ have
negative orders at each valuation at infinity of $f_Y$ and hence our
reductions are simply valid.
Thus we get the global Milnor number
$$\mu(f;A) = \sum_{a \in k} \mu(f-a) = 1+r_0 - \sum_{1}^h (n_i-1)r_i 
=-C(\Gamma(f)).$$

\begin{proposicion} 
If $f\in A$ is a curve with one place at infinity, then we have
$$\mu(f;A)=C(f)+2P_g(f) = \mu(f)+\overline{\chi}(f)+2P_g(f) 
=\nu(f)+\delta(f)+\overline{\chi}(f)+2P_g(f) $$

\end{proposicion}
\begin{demostracion}{.}
This follows by using already obtained results from Lemmas in Section 3 as well
as conclusions of Sections 4.5 and 4.6.

\end{demostracion}

\begin{corolario} 
If $f$ is a curve with one place at infinity such that $\mu(f;A)=\nu(f)$
then we have that $f$ is a rational curve which is unibranch and locally 
quasihomogeneous at all its points. Further all translates of $f$ are
nonsingular. Moreover, 
$Z(f)=l(\Omega(\overline{R},k)/\Omega(R,k))=C(f)/2$.
\end{corolario}

\begin{demostracion}{.}
It follows from the Proposition 4.1 that 
$\delta(f)+\overline{\chi}(f)+2P_g(f)=0$ and since all the quantities
are nonnegative, each is zero. This proves the first assertion. 
It follows that
$\mu(f;A)=C(f) = \mu(f)$. 
Thus $\mu(f-a)=0$ for all $a \neq
0$ giving that all translates of $f$ are nonsingular. 
The last result follows from $C(f)=\nu(f)$ using Proposition 2.4.
\end{demostracion}

\subsection{ Connection with the NP (Newton-Puiseux) expansion.}
All the above quantities can either inductively be described and
calculated from the known valuation, or they can be deduced from the
NP expansion which defines the valuation at infinity.

The NP expansion is a power series parametrization:
$$x=\tau^{-n}~,~ y=\eta(\tau)=\nz \tau^{-m}+\cdots \in k((\tau)).$$
The valuation $v$ can then be described as
$v(p(x,y))=\ord_\tau(p(\tau^{-n},\eta(\tau))$.

We refer the reader to look up any of the following sources: 
\cite{A1},\cite{SS},\cite{S},\cite{AM1},\cite{AM2}, for the definition of the
associated characteristic sequences and calculations.

For future use, we prove:
\begin{lema} 
Suppose that $f\in A$ has one place at infinity and an NP expansion
$$x=\tau^{-n}~,~ y=\eta(\tau)=\nz \tau^{-m}+\cdots \in k((\tau)).$$
Then there are no terms of the form $\nz \tau^{ns}$ in the expansion of
$\eta(\tau)$ for $s >0$.

\end{lema}

\begin{demostracion}{.}
Write $\eta(\tau) = \sum_{-m}^{\infty} a_i \tau^i$.
We know that 
$$f=\prod_{\zeta^n =1} (Y-\eta(\zeta \tau)).$$
and the coefficient of $Y^{n-1}$ is then easily seen to be the trace
$\displaystyle{-\sum_{\zeta^n=1} \eta(\zeta \tau)}$
where we substitute $X$ for $\tau^{-n}$. Suppose if possible,
$s>0$ and $a_{ns}\neq 0$. Then
we get that the term $Y^{n-1}X^{-s}$ has the nonzero 
coefficient 
$\displaystyle{-a_{ns}\sum_{\zeta^n=1}(\zeta^{ns})=-a_{ns}(n)}$. Since the 
coefficient of $Y^{n-1}$ must be a polynomial in $X$, we get a contradiction, 
thereby proving the result.

\end{demostracion}

\section{Rational one place curves}

Let the notations be as in the Section above. Assume that $f$ defines 
a rational curve. It follows by a modification of the L\"{u}roth theorem
that it is a ``polynomial curve'' which means that there are two polynomials 
$X(T),Y(T)$ in a new indeterminate $T$, such that $f(X(T),Y(T))=0$ and
at least one of $X(T),Y(T)$ is a non constant.
 
As before, let us take the image $\phi_f(T)=t$ and write
$x(t)=\phi_f(X(T)), y(t)=\phi_f(Y(T))$, so that $f(x(t),y(t))=0$.
By choosing minimal $t$-degrees for $x,y$,
we may assume $k(f)=k(t)$ and hence  $\overline{R}=k[t]$. 
It follows that $P_g(f)=0$ and hence 
$C(\Gamma(f)) = -C(f)$.
Also we have 
$$l(\overline{R}/R)=
l(R/\mathfrak{C})= C(f)/2 = -C(\Gamma(f))/2.$$

The (unique) valuation $v$ at infinity of $f$ is now seen to be the
$1/t$-adic valuation which can be simply described as
$v(\frac{a(t)}{b(t)})= \deg_t(b(t))-\deg_t(a(t))$.
The order  $v(a(t)db(t))$ can now be seen as $v(a(t))+v(b(t))-1$
as long as $a(t)d(b(t)) \neq 0$, i.e. $0\neq a(t)\in k(t)$ and 
$b(t)\in k(t)\setminus k$.

We easily see that the set 
$$\{v(h) ~|~ 0 \neq h \in \overline{R}=k[t]\} = -N $$
the set of all negative natural numbers. The semigroup of the curve is a
subset of this with exactly $(-1/2)C(\Gamma(f))=(1/2)C(f)$ missing
values.

We now prove:
\begin{lema}
With the notation as above, let
$\Gamma^{*}(f)=\{ v(p) ~|~  p \in \Omega(R,k) \}$ the
set of values of elements in $\Omega(R,k)$ and 
$\Gamma'(f) = \{ v(u)-1 ~|~ u \in R \setminus k\}$, the set of 
values of exact differentials $du\in \Omega(R,k)$.
Further, set:
$M = -N \setminus \Gamma^{*}(f)$.
We have the following:

\begin{enumerate}
\item $\Gamma'(f) \subseteq \Gamma^{*}(f)$.

\item 
$Z(f)=$ the number of elements of $M$.
\item
The cardinality of $-N \setminus \Gamma'(f)$ is clearly $-C(\Gamma(f))/2$ and hence
$Z(f) \leq -C(\Gamma(f))/2$ and the inequality is strict if and only if 
$\Gamma^{*}(f) $ is strictly bigger than $\Gamma'(f)$ or equivalently
there is some element $p \in \Omega(R,k)$ which is not exact.
\end{enumerate}

\end{lema}

\begin{demostracion}{.}
The first statement is obvious.

Given any element $a\in \Omega(\overline{R},k)=k[t]dt$ we claim that
\begin{itemize}
\item $a\in \Omega(R,k)$ or,
\item there is an element $b\in \Omega(R,k)$ such that $v(a-b)\in M$. 
\end{itemize}

If $a \not \in \Omega(R,k)$, then we can choose $b\in \Omega(R,k)$
which maximizes the value of $v(a-b)$, since all values in question are
bounded above by $0$. If $v(a-b) \in M$ then we are done.
Otherwise, there exists $b^*\in \Omega(R,k)$ such that
$v(a-b) = v(b^*)$. Then $v(a-b-\lambda b^* )>v(a-b)$ for some 
$\lambda\in k$, contradicting the choice of $b$. 

Hence $v(a-b) \in M$ as claimed. 

This shows that we can make a basis of
$\Omega(\overline{R},k)/\Omega(R,k)$ consisting of elements with
distinct values in $M$, proving our first second assertion.

The last assertion easily follow.
 
\end{demostracion}

We now summarize the results.

\begin{proposicion} Let $f\in A$ be a  curve with one place at infinity.
Then $\mu(f;A)=\nu(f)$ if and only if
$f$ is a rational curve which is locally unibranch at all its points
and all differentials in $\Omega(R,k)$ are exact.

\end{proposicion}
\begin{demostracion}{.}
We have already proved the only if part in view of Corollary 4.2 
except for the claim of exactness of differentials. This claim follows
from Lemma 5.1, since under our hypothesis we know that
$-C(\Gamma(f))=C(f)$.

The if part easily follows from the Proposition 4.1, since our
hypothesis implies that the quantities
$P_g(f),\overline{\chi}(f),\delta(f)$ are all zero.
\end{demostracion}

{\bf We now define} the following terms.
Let $z_1,z_2$ be non constant elements of the function field $k(f)$
of an irreducible curve and let $w$ be any valuation of $k(f)$ over $k$.
The differential 
$\beta(z_1,z_2;w)= w(z_1)z_1dz_2 - w(z_2)z_2dz_1$
is defined to be {\bf the basic differential associated with the pair 
$(z_1,z_2)$ corresponding to the valuation $w$}. 

We also define  {\bf the gap of the basic differential}
$\beta(z_1,z_2;w)$ to be

$$w(\beta(z_1,z_2;w)) - w(z_1)-w(z_2)+1.$$
In case the basic differential is zero the gap is defined to be
$\infty$.

If we take $x,y$ to be the generators of the coordinate ring, then we will
show in the next Lemma 5.3 that the gap of $\beta(x,y;v)$ is completely
determined by the NP expansion associated with $x,y$. We shall say that
the pair $(x,y)$ has {\bf a maximal gap }if for all automorphisms $\sigma$ of 
$A$ such that $v(\sigma(x))=v(x),~ v(\sigma(y))=v(y)$, the gap 
of the basic differential $\beta(\sigma(x),\sigma(y);v)$ is at most 
equal to the gap of $\beta(x,y;v)$.
Note that an infinite gap is clearly maximal.

We now prove an auxiliary Lemma.
\begin{lema}
Suppose that $f\in A$ has one place at infinity and an NP expansion
$$x=\tau^{-n},y=\eta(\tau) = \tau^{-m} + \nz \tau^{-m+q} + \mbox{ higher
terms.}$$
Let $v$ denote the associated valuation at infinity.

Assume that $(n,m)$ is non principal, i.e.
neither divides the other.

Then the gap of $\beta(x,y;v)$ is exactly $q$.

Suppose that the pair $(x,y)$ has a maximal gap as defined above.

Then we have the following:

\begin{enumerate}
\item $-m+q$ is not a multiple of $n$.
\item $-n+q$ is not a multiple of $m$.
\item An equation $-m-n+q = -an -bm $ does not hold for non negative
integers $a,b$.
\end{enumerate}

\end{lema}

\begin{demostracion}{.}
The calculation of the gap of the basic differential is
straightforward substitution. It follows that the value of $q$ does not
change by interchanging the role of $x,y$.

The proofs are by contradiction.
We show that when the desired conditions on $q$ fail, we can make an
automorphism $\sigma$ on $x,y$ so that the gap
of the basic differential $\beta(\sigma(x),\sigma(y);v)$ is bigger than
$q$ but $v(\sigma(x))=v(x)$ and $v(\sigma(y))=v(y)$. This would give a 
contradiction.

Consider the first case and suppose is possible $-m+q=-ln$.
First, we claim that $l\geq 0 $, for otherwise, we have a contradiction by
Lemma 4.3. 
Now consider new variables $X'=X,Y'=Y+uX^l$ where $u \in k$.
Let $x',y'$ be their respective images modulo $f$ and note that:
$$v(x)=v(x')~,~v(y)=v(y')~,~
\beta(x',y';v) = \beta(x,y,v)+(v(x)l-v(y))u x^ldx.$$ 

It is easy to verify that the leading term of $\beta(x',y';v)$ is
$(\nz + \nz u)\tau^{-ln-n-1}d\tau$  and so by choosing $u\in k$ to kill this
term we can make its order bigger than $-ln-n-1=-m-n+q-1$. It follows
that the gap gets bigger than $q$. This is a contradiction to the
assumed maximality of $q$.

The second case has a similar argument.

Now in the third case, we already know that if $-m-n+q = -an -bm $
then both $a,b$ are non zero. it follows that the right hand side is
less than or equal to $-n-m$ giving $q \leq 0$. Since $q>0$ we have a
contradiction!

\end{demostracion}

\begin{corolario}
Suppose that $f\in A$ has one place at infinity and an NP expansion
$$x=\tau^{-n},y=\eta(\tau) = \tau^{-m} + \nz \tau^{-m+q} + \mbox{ higher
terms.}$$

Further assume that $(n,m)$ is non principal.

Let $v$ denote the associated valuation at infinity.
Then there is an automorphism $\sigma$ of $A$ such that
$v(\sigma(x))=v(x),~ v(\sigma(y))=v(y)$ and the pair 
$(\sigma(x),\sigma(y))$ has a maximal gap as defined above.

Further, the automorphism $\sigma$ is ``very elementary'' which means 
that it is either of the form
$\sigma(X)=\nz X, \sigma(Y) = \nz Y+l(X)$  or of the form
$\sigma(X)=\nz X+l(Y),\sigma(Y)=\nz Y$ with $l(X)\in k[X]$.

\end{corolario}

\begin{demostracion}{.} From known structure of automorphisms of curves with one place at
infinity as in \cite{AM3} and \cite{ASI}, we see that automorphisms
which preserve that $v$-values of $x,y$ are necessarily very elementary.

Considering how the equation $f(X,Y)$ is related to the NP expansion, we
see that either $q=\infty$ and we have the desired maximal gap or $q$ is
bounded above by $mn-1$. 

The proof of Lemma 5.3 shows that the gap of the basic differential
$\beta(x,y;v)$ is precisely $q$ and since this is bounded above 
there is a clear maximum value for this $q$ under all possible very 
elementary automorphisms.

We simply choose an automorphism which gives such a maximal value.
\end{demostracion}

\noindent We now prove:
\begin{teorema}{\bf Generalized Theorem of Zariski.}  
Let $f\in A$ be a curve with one place at infinity.
Then $\mu(f;A)=\nu(f)$ if and only if 
there is an automorphism $\sigma : A \longrightarrow A$ such that 
$f=\sigma(Y)^b-\sigma(X)^a$ for non negative integers $a,b$.
Moreover, in view of irreducibility of $f$, we must have $\mbox{GCD}(a,b)=1$.

\end{teorema}
\begin{demostracion}{.} It is easy to check that when $f$ has the indicated
form $\sigma(Y)^b-\sigma(X)^a$, we have
$$\mu(f;A) = l(A/(\sigma(X)^{a-1},\sigma(Y)^{b-1})A = (a-1)(b-1)=\nu(f).$$
All the remaining properties of $f$ are easy to verify.

Now for the only if part, we note that Proposition 5.2 already gives
that $f$ is a rational unibranch curve with $\delta(f)=0$ and that all
differentials in $\Omega(R,k)$ are exact. We use this to prove the
existence of the automorphism.

We can write after an adjustment of coordinates:
$$f(X,Y) =  Y^n -X^m + \sum f_{ij}X^iY^j $$
where $n >m$ and the terms of the summation are naturally restricted by the
condition that $f$ has one place at infinity.

If $m=0$, then we have $y\in k$ and after a translation in $Y$ we get
$f=Y$. This is the case when $a=1,b=0$.

Now assume that $n>m>0$ such that $(n,m)$ is non principal, i.e.
neither divides the other.
This assumption is standard and ensured by performing elementary
automorphisms on $X,Y$ if one of $n,m$ divides the other.

We shall now assume the notations and details about the Abhyankar-Moh theory
in Section 4 without further explanation.

In particular, $r_0=-n,~ r_1=-m,~ d_1 =n,~ d_2 = \mbox{GCD}(m,n)>1$ and by
our arrangement $n>m>d_2$.

Let an NP expansion be written as:
$$x=\tau^{-n} ~,~ y=\eta(\tau) = \tau^{-m}  + \mbox{ higher terms} $$
and let $v$ denote the corresponding valuation at infinity.

If there are no higher terms, then we have $y^n=x^m$ and clearly
$f=Y^n-X^m$. By irreducibility of $f$ we get $\mbox{GCD}(m,n)=1$
finishing the proof.

So, now assume that:
$$x=\tau^{-n} ~,~ y=\eta(\tau) = \tau^{-m} +\nz \tau^{-m+q} + 
\mbox{ higher terms} $$

As in Corollary 5.4, we assume that the gap $q$ of the basic differential 
$\beta(x,y;v)$ can be chosen to be maximal
so that we can assume that the conditions established in Lemma 5.3 hold.

By exactness of all differentials in $\Omega(R,k)$, we deduce that 
$-m-n+q \in \Gamma(f)$.

We have two possible cases and we show that both fail.
\begin{enumerate}
\item{Case 1:~~} $-m+q$ is the second characteristic exponent, i.e.
$q=q_2$. Since $r_2 = n_1r_1 + q_2$,  we have
$$-m-n+q = r_0 + r_1+r_2 - n_1 r_1 = r_0 + r_1+r_2 - n_0 r_0 =
(1-n_0)r_0+r_1+r_2.$$

This is an admissible expression and hence is in $\Gamma(f)$
if and only if $1-n_0 \geq 0$. Further, this is
possible only if $n_0 = 1$, i.e. $r_0=-n$ divides $r_1=-m$. This is a
contradiction.

\item{Case 2:~~} $-m+q$ is not a characteristic term, so $-m+q$ and hence
$-m-n+q$ is divisible by $d_2$ and we have
$-m-n+q = a_0 r_0 + a_1 r_1$ where $0 \leq a_1 < n_1$ and $0 \leq a_0$.

Thus $-m-n+q = -a_0n - a_1 m$ in contradiction with condition 3 of Lemma
5.3.

\end{enumerate}

\end{demostracion}

\section{Link with Lin-Zaidenberg Theorem}

\medskip

\noindent Let $f\in A$ be an irreducible polynomial.
In \cite{LZ},
Lin and Zaidenberg proved the following: 

\begin{teorema}{\bf Lin-Zaidenberg} If $f$ is a rational curve with one place at infinity
which is locally unibranch at all its points then 
there is an automorphism $\sigma : A \longrightarrow A$ such that 
$f=\sigma(Y)^b-\sigma(X)^a$ for non negative integers $a,b$.
Moreover, in view of irreducibility of $f$, we must have $\mbox{GCD}(a,b)=1$.

Moreover, the hypothesis of the theorem can be alternatively stated as:
$f\in A$ is a curve with one place at infinity such that
$\mu(f;A)=\mu(f)$.

\end{teorema}

Their results generalized the celebrated Abhyankar-Moh Epimorphism Theorem
which can be stated as follows:

\begin{teorema}{\bf Abhyankar-Moh} If $f$ is a rational nonsingular curve 
with one place at infinity then 
there is an automorphism $\sigma : A \longrightarrow A$ such that 
$f=\sigma(Y)$. 
\end{teorema}

\begin{nota}
{\rm 
Note that the hypothesis of the Lin-Zaidenberg Theorem implies in our
notation:
$P_g(f)=0=\overline{\chi}(f)$. Then, in view of sections 4.5, 4.6 we have
$$\mu(f;A) = C(f) = \mu(f)=\nu(f)+\delta(f).$$

In turn, the alternate hypothesis
$\mu(f;A) = \mu(f)$ is easily seen to imply all their conditions as in
our proof of Corollary 4.2. 

The Epimorphism Theorem is a special case in view of the facts that
$\mu(f)=0$ by the non singularity of the curve while the formula in 4.6
gives:
$\mu(f;A)=-C(\Gamma(f))=0$. Thus it is now a special case of the
Lin-Zaidenberg Theorem. We note that, since the proofs of the Lin-Zaidenberg 
theorem do use the Epimorphism Theorem, they do not provide an
alternative proof for it.

It follows from the conclusion of the Lin-Zaidenberg theorem that 
further,  $\delta(f)=0$. All known 
proofs of the Lin-Zaidenberg theorem are rather involved, using
respectively complex analysis (Lin and Zaidenberg \cite{LZ}), 
topology (Neumann and Rudolph \cite{NR}),  
and finally algebraic geometry using surface theory 
(Gurjar and Miyanishi \cite{GM}). 

Our original aim was to provide a much simpler proof of this important
theorem using the techniques of the Abhyankar-Moh theory.

Our theorem is a weaker version of the Lin-Zaidenberg theorem since we
additionally need to assume that $\delta(f)=0$ and this is clearly the
condition that $J_f = Jac_f$ or $f \in J_f$. The Lin-Zaidenberg
hypothesis easily implies that $f \in \sqrt{J_f}$ and the Lin-Zaidenberg
theorem implies that $f \in J_f$. We have not yet
succeeded in establishing this independently using our simpler 
techniques. This condition can be also described as
the condition of local quasihomogeneity explained
earlier in Remark 3.5. 
}
\end{nota}

\section{Estimating independent non exact differentials}
\subsection{Preamble}
Let $f\in A$ be a rational curve with one place at infinity.
We know from Proposition 4.1 and Lemma 3.4 that for such a curve
$C(f) = -C(\Gamma(f))$ and 
$C(f)/2-Z(f) = \delta(f)+\overline{\chi}(f)$ measures how far 
the curve fails to be quasihomogeneous.
Also the quantity $C(f)/2-Z(f)$ is seen to be the number of independent
non exact differentials by calculations of Lemma 5.1.
\footnote {In the notation of the Lemma 5.1 this is the cardinality of
the difference $\Gamma^{*}(f) \setminus \Gamma^{'}(f)$.}

In this section we establish some estimates of this quantity.

In the above notation we assume:
$f\in A$ is a rational curve with one place at infinity with  NP expansion
$$
X=\tau^{-n}~,~
Y=\tau^{-m}+\nz \tau^{-m+q}+\cdots
$$
and let $v$ as usual denote the resulting valuation. Note that we are
implicitly assuming that $q \neq \infty$, in other words, our curve
is not globally quasihomogeneous.

We further assume that $n>m>\mbox{GCD}(n,m)$ and the gap $q$ is chosen
to be as large as possible for the given coprime degrees. 

We already know that the basic differential
$\Delta = \beta(x,y;v)$ associated with $x,y$ is not
exact in view of our proof of Theorem 5.5.

We investigate the set of non exact differentials which are multiples of
this $\Delta$.

By the two cases analyzed in the proof of our Theorem 5.5, we see that
we have:
\begin{enumerate}
\item {Case 1:~~} $-m+q=-m+q_2$ the second characteristic number and $-m-n+q =
(1-n_0)r_0+r_1+r_2 \not \in \Gamma(f)$ since $n_0 >1$.
\item{Case 2:~~} $-m+q$ is not a characteristic term and hence
$-m-n+q = a_0 r_0 + a_1r_1 $ where $0 \leq a_1 <n_1$ and $a_0 <0$.
\end{enumerate}

\subsection{Number of distinct values of inexact differentials}
We now prove 

\begin{lema} 
Assume the setup of the preamble above and assume that we have case 1.
Let $w_b$ be any monomial $\prod g_i^{b_i}$ where $b=(b_0,\cdots,b_h)$ is
an admissible sequence satisfying 
$$0 \leq b_0 <n_0-1, 0\leq b_1 <n_1-1, 0\leq b_2<n_2-1,0\leq b_i <n_i
\mbox{ for }3 \leq i \leq h.$$

Then the differentials $\{ w_b\Delta\}$ are inexact members of 
$\Omega(R,k)$ with distinct values. In particular, the number of these
is at least $(n_0-1)(n_1-1)(n_2-1)\prod_3^h n_i \geq 2^{h-1}$.
\end{lema}

\begin{demostracion}{.}
The differentials are inexact because their values augmented by $1$ are
not in $\Gamma(f)$ as evident from their standard representation and it 
also gives the distinctness of the values. The estimate on the count is
simply by counting the possible values of $b$'s. For the last estimate 
note that each $n_i$ is at least $2$. Moreover the
first two $n_0,n_1$ are coprime to each other and hence at least one of
them is $3$ or bigger. This gives that at least $h-1$ of the $h+1$
factors are at least $2$ and finishes the proof.
\end{demostracion}

\begin{lema} 
Assume the setup of the preamble above and assume that we have case 2.
Let $w_b$ be any monomial $\prod g_i^{b_i}$ where $b=(b_0,\cdots,b_h)$ is
an admissible sequence satisfying 
$$0 \leq b_0 <-a_0, 0\leq b_1 <n_1-a_1, 0\leq b_i <n_i
\mbox{ for }2 \leq i \leq h.$$

Then the differentials $\{ w_b\Delta\}$ are inexact members of 
$\Omega(R,k)$ with distinct values. In particular, the number of these
is at least $(-a_0)(n_1-a_1)\prod_2^h n_i \geq 2^{h-1}$.
\end{lema}

\begin{demostracion}{.}
Everything except the last statement follows exactly as in Lemma 7.1.

The last estimate 
follows from the fact that each $n_i$ is at least $2$ and we have 
$h-1$ terms $n_2,\cdots,n_h$ as factors of our estimate.
\end{demostracion}

We now investigate the case when the number of values of inexact differentials
is exactly $1$. We prove:

\begin{proposicion} 
Suppose that $f\in A$ is a curve with one place at infinity such that 
$\mu(f;A)=\nu(f)+1$. Then we have the following:
\begin{enumerate}
\item $f$ is a rational curve with $\delta(f)=0$ and
$\overline{\chi}(f)=1$.
\item There is exactly one  value of an inexact differential in $\Omega(R,k)$.
\item All such curves have $h=1$ or equivalently $n,m$ are coprime.
Moreover, we have exactly one of the following three situations:
\begin{itemize}
\item{\bf Situation 1:} $m=2,n=2p+1$ for some $p=1,2,\cdots$.
\item{\bf Situation 2:} $m=3,n=4$.
\item{\bf Situation 3:} $m=3,n=5$.
\end{itemize}

Each of the these situations will be explicitly described below.
\end{enumerate}

\end{proposicion}
\begin{demostracion}{.} From the given condition we deduce from Proposition 4.1 that 
$$\delta(f)+\overline{\chi}(f)+2P_g(f)=1.$$
It follows that $P_g=0$, i.e. the curve is rational. If
$\overline{\chi}(f)=0$, then we could apply the Lin-Zaidenberg Theorem 6.1 and
deduce that $\delta(f)=0$ giving a contradiction. This gives the first
part. 

The second part follows from Lemma 3.4. 

>From the lower bound $2^{h-1}$  in Lemma 7.1, 7.2 we deduce that $h=1$. 
Further, this means that we 
can apply the setup of Lemma 7.2 and in the notation of the preamble
we get that $$\mbox{GCD}(m,n)=1, a_0=-1, a_1=n_1-1=n-1 \mbox{ and }
-m -n +q = (-1)(-n)+(n-1)(-m).$$
By adding $4$ to both sides and rearranging, the last equation reduces to
$$nm-2n-2m +4 =(n-2)(m-2)= 4-q  $$

Since $n>m$ by the preamble, we see that $m=2,3$, since otherwise
$4<(n-2)(m-2)=4-q<4$ is a contradiction.

If $m=2$ then we get $q=4$ and $n=2p+1$ for some $p=1,2,\cdots$ since it is
coprime with $m$.

If $m=3$, then $n=4,q=2$ and $n=5,q=1$ give the only possibilities,
since $n>5$ leads to $4 \leq (n-2)(m-2) = 4-q <4$ a contradiction again!

We shall now give the explicit description of each of these cases below in
separate Lemmas.
\end{demostracion}

\begin{lema} 
Assume that we have Situation 1 as described in Proposition 7.3.
Then by an automorphism we can arrange that $f=g(Y)-X^2$ where
$g(Y)\in k[Y]$ is monic of degree $2p+1$ and has exactly two
distinct linear factors. In turn, every such curve satisfies the
hypothesis of Proposition 7.3.
\end{lema}

\begin{demostracion}{.}
Since $m=2$ we know that $f$ is of degree 2 in $X$ and has a non zero
constant coefficient for $X^2$.
Thus the form of $f$ can easily be arranged by translating $X$ by a
polynomial in $Y$ to kill the $X$-term.
Now let $s$ be the number of distinct linear factors of $g(Y)$. We see:
$$Jac_f = (f,f_Y,f_X)A = (f,g'(Y),-2X)A = (X,g(Y), g'(Y))A.$$
It is not hard to see that the length of the last ideal 
is $\nu(f)=2p+1-s$ and since this
is required to be $\mu(f;A)-1 = 2p-1$ we get that $s=2$. In turn, this
also shows that all such curves satisfy the hypothesis of Proposition
7.3.
\end{demostracion}

\begin{lema} 
Assume that we have Situation 2 as described in Proposition 7.3.
Then we can arrange the parametric form of the curve as
\begin{itemize}
\item  either $x=t^4+at^2,y=t^3-3t^2, a=0,-8,-9$ 
in case there is at least one singular
(cuspidal) branch at finite distance 
\item or $x=(t^2-1)(t^2+a),y=(t^2-1)(t+b)$, where $a=-1$ and $ b \neq \pm 1$ in
case there is no singular branch. 
\end{itemize}

\end{lema}

\begin{demostracion}{.}
The proof is done by using Maple and we will only outline the strategy.
\footnote{ A Maple file with actual calculations will be provided if needed.}

First, we assume that the singular branch is at $t=0$ and then we can
easily assume that either $y=t^3$ or 
$y=t^3+\nz t$. 

In the first case, we have $x=t^4+\nz t^2$ and a detailed calculation
leads to impossible subcases.

In the second case, by a suitable scaling we may arrange $y=t^3-3t^2$.
 
By singularity of the branch,
$x$ has order at least $2$ and we can remove its cubic term by adding a
multiple of $y$.

We next calculate the Taylor Resultant (TRES) introduced by Abhyankar and
described thus:

Given a parametric curve $x=u(t),y=v(t)$ the branches at the singular
points are determined by calculating the 
$$TRES(u(t),v(t))=
\mbox{Resultant}(\frac{u(t)-u(s)}{t-s},\frac{v(t)-v(s)}{t-s},s).$$

This resultant has precisely the degree equal to $C(f)$ and the 
multiplicities of its factors give the
exact order of the conductor ideal at various branches i.e.  valuations.
Indeed, for a rational one place curve, it gives the exact generator of
the conductor ideal $\mathfrak{C}_f$ in the ring $k[t]$.

Under our hypothesis of $\overline{\chi}(f)=1$, we see that this
 TRES has degree $6$ and has at most two simple roots and thus a 
maximum of $4$ distinct roots. Systematic use of the discriminant lets
us determine the conditions on the coefficients $(a,b)$ leading to the
announced values.

For the case when we don't have a singular (cuspidal) branch, we assume 
that two branches (which can be arranged without loss of generality to
be ) $t=\pm 1$ are centered at a singular point arranged to be the
origin. We then get the announced parametrization after killing out the cubic
term by adding a multiple of $y$ to $x$. We now note that there cannot
be another singular point, for it necessarily will be unibranch and will 
give a singular branch reducing to the first case. Also, at the origin,
we must have just the two branches and three nodes in successive
neighborhoods. Calculations with usual quadratic transformations give
the announced result.

\end{demostracion}

\begin{lema} 
Assume that we have Situation 3 as described in Proposition 7.3.
Then we can arrange the parametric form of the curve as
\begin{itemize}
\item  either there is at least one singular
(cuspidal) branch at finite distance 
and $$x=t^5+at^4+bt^2,y=(t^3-3t^2)$$ where the values 
of $(a,b)$ come from the finite set: 
$$\{(-15,0),(-3,0)(-5/2,0),(-3/2,0),(-3,4),(-6,27),(-7,36),(-15/2,40)\};$$
\item or there is at least one singular
(cuspidal) branch at finite distance 
and 
$$x=t^5+as^4+bs^2,q=s^3 \mbox{ with } b=0,a\neq 0;$$
\item or $$x=(t^2-1)(t^3 \pm \sqrt{2} t^2+1),y=(t^2-1)(t+1 \pm \sqrt{2} )$$ 
in case there is no singular branch. 
\end{itemize}
\end{lema}

\begin{demostracion}{.}
The arguments here are similar to Lemma 7.5, except we get an extra
variable in each case.

In case of at least one singular branch, we argue as before and 
both the subcases now lead to 
indicated solutions.

In case there is no singular branch, we again make a similar argument 
and end up with only a finite set of solutions as indicated. Note that
in this case, we need four successive nodes at the origin (instead of the
three nodes in Lemma 7.5) and hence the
loss of a free parameter is expected, so we end up with only a finite set
of solutions.

The reduction in the first case is identical and the argument in the
second case that the origin must be the only singular point still holds,
except we need to find four nodes in successive neighborhoods.

\end{demostracion}

\begin{nota}{\rm It seems difficult to generalize the characterization 
above to one place rational curves with small $\mu(f;A)-\nu(f)$, 
for example, if $\mu(f;A)-\nu(f)=2$, then Lemma 7.1. and Lemma 7.2. 
show that $h=1$ or $2$. Every case contains many subcases. This makes
the (already technical) work hard to realize.}
\end{nota}


\newpage 
\begin{thebibliography}{99}

\bibitem[A1]{A1}
Abhyankar, S. S.
\newblock {\it ``Expansion Techniques in Algebraic Geometry'',}
\newblock Tata Institute of Fundamental Research, 1977.

\bibitem[AM1]{AM1}
Abhyankar, S.S., and Moh, T.T., 
\newblock {\it Newton-Puiseux expansion and Tschirnhausen
transformation I,}
\newblock J. Reine Angew. Math. 260 (1973), 47-83.

\bibitem[AM2]{AM2}
Abhyankar, S.S., and Moh, T.T., 
\newblock {\it Newton-Puiseux expansion and Tschirnhausen
transformation II,}
\newblock J. Reine Angew. Math. 261 (1973), 29-54.

\bibitem[AM3]{AM3}
Abhyankar, S.S., and Moh, T.T., 
\newblock {\it Embeddings of the line in the plane,}
\newblock J. Reine Angew. Math. 276 (1975),148-166.

\bibitem[ASI]{ASI}
Abhyankar, S.S., and Singh,B., 
\newblock{\it Embeddings of certain curves in the affine plane,}
\newblock Amer. J. Math. 100(1978), 99-175.

\bibitem[AS]{AS}
Abhyankar, S.S., and Sathaye,A., 
\newblock{\it Geometric Theory of Algebraic Space Curves,}
\newblock Lecture Notes in Mathematics, 423(1974).


\bibitem[AS2]{AS2}
Abhyankar, S.S., and Sathaye,A., 
\newblock{\it Uniqueness of plane embeddings of special curves,}
\newblock Proceedings of the AMS, 4 (1996), 1061-1069.


\bibitem[BR]{BR}
Brian\c{c}on, J.
\newblock{\it Un polyn\^ome appartien-il \`a l'id\'eal
de ses d\'eriv\'ees partielles?}
\newblock Preprint n$^{o}$632, Nov 2001, Universit\'e de Nice-Sophia Antipolis.

\bibitem[GM]{GM}
Gurjar, R. V., and Miyanishi M., 
\newblock{\it On Contractible Curves in the Complex Affine Plane}
\newblock Tohuku Math. J(2) 1996, n$^{o}$3, 459-469.

\bibitem[LZ]{LZ}
Lin, V. I., and Zaidenberg,M., 
\newblock{\it An irreducible simply connected algebraic curve in ${\bf C}^2$
is equivalent to a quasihomogeneous curve}
\newblock Dokl. Akad. Nauk. SSSR 271(1983), n$^{o}$5, 1048-1052.

\bibitem[NR]{NR}
Neumann, W. and Rudolph, L.
\newblock{\it Unfoldings in knot theory (and Corrigendum)}
\newblock  Math. Ann. 278(1987) and 282(1988), 409-439 and 349-351.

\bibitem[S]{S}
Sathaye, A., 
\newblock {\it Generalized Newton-Puiseux expansion and Abhyankar-Moh
semigroup theorem,}
\newblock Invent. math. 74(1983), 149-157.

\bibitem[SS]{SS}
Sathaye, A. and Stenerson Jon,
\newblock {\it On Plane Polynomial Curves.}
\newblock {\it  Algebraic Geometry and Applications, (1994), 121-142.}

\bibitem[SA]{SA}
Saito, K., 
\newblock {\it Quasihomogen isoliere Singularit\"{a}ten von
Hyperfl\"{a}chen,}
\newblock Invent. math. 14(1971), 123-142.

\bibitem[Z1]{Z1}
Zariski, O., 
\newblock {\it Characterization of Plane Algebroid Curves whose module
of Differentials Has Maximum Torsion,}
\newblock Collected Works of Oscar Zariski, Vol. III, 475-480.

\end{thebibliography}
\end{document}